\documentclass[a4paper,11pt]{article}
\usepackage{amsfonts}
\usepackage{amsmath}
\usepackage{hyperref}

\input{tcilatex}
\begin{document}

\title{Estimating the error term in the Trapezium Rule using a Runge-Kutta
method}
\author{J. S. C. Prentice \\
Senior Research Officer\\
Mathsophical Ltd.\\
Johannesburg, South Africa\\
Email: jpmsro@mathsophical.com}
\maketitle

\begin{abstract}
We show how the error term for the Trapezium Rule can be estimated, by
solving an initial value problem using a Runge-Kutta method. The error term
can then be added to the Trapezium approximation, yielding a much more
accurate result. We also show how the risk of singularities in the relevant
initial value problem can be mitigated.
\end{abstract}

\section{Introduction}

Students of applied mathematics, and numerical analysis in particular, have
no doubt encountered the Trapezium Rule \cite{burden}\cite{kincaid} and
Runge-Kutta (RK) methods \cite{butcher}. They will very likely have met the
error term for the Trapezium Rule, and they may even know how to estimate an
integral using a RK method. However, they may not know that the Trapezium
Rule's error term can be estimated by solving a suitable initial value
problem using, of course, a RK method. Once the error term is known, the
Trapezium approximation can be considerably improved. This interesting and
potentially useful device is the topic of this paper.

\section{Theory}

We have the equation%
\begin{eqnarray}
\int\limits_{a}^{x}f\left( x\right) dx &=&\frac{\left( x-a\right) }{2}\left(
f\left( a\right) +f\left( x\right) \right) -\frac{\left( x-a\right) ^{3}}{12}%
f^{\prime \prime }\left( \xi \left( x\right) \right)  \label{problem} \\
\xi &\in &\left( a,x\right) .  \notag
\end{eqnarray}%
Here, $f:%
\mathbb{R}
\rightarrow 
\mathbb{R}
$ is assumed to be Riemann integrable, the upper limit in the integral is
variable, the first term on the RHS is the Trapezium Rule written in terms
of $x,$ and the second term on the RHS is the error in the Trapezium Rule,
wherein we assume that $\xi $ is a function of $x$.

Differentiating wrt $x$ gives%
\begin{align*}
f\left( x\right) =\text{ }& \frac{1}{2}\left( f\left( a\right) +f\left(
x\right) \right) +\frac{\left( x-a\right) f^{\prime }\left( x\right) }{2} \\
& -\frac{3\left( x-a\right) ^{2}}{12}f^{\prime \prime }\left( \xi \right) -%
\frac{\left( x-a\right) ^{3}}{12}\frac{df^{\prime \prime }\left( \xi \right) 
}{d\xi }\frac{d\xi }{dx}
\end{align*}%
which yields the differential equation (DE)%
\begin{equation}
\frac{d\xi }{dx}=\frac{-18f\left( x\right) +6f\left( a\right) +6\left(
x-a\right) f^{\prime }\left( x\right) -3\left( x-a\right) ^{2}f^{\prime
\prime }\left( \xi \right) }{\left( x-a\right) ^{3}f^{\prime \prime \prime
}\left( \xi \right) }.  \label{DE}
\end{equation}%
If an initial value $\xi \left( x_{0}\right) $ is provided, this DE can be
solved using a Runge-Kutta method to yield $\xi \left( x\right) .$ This
allows the error term to be computed, which can then be added to the
Trapezium term to obtain a more accurate approximation for the integral.

Note that we must necessarily assume that $f$ and its first three
derivatives all exist on $\left[ a,x\right] .$

\section{Example}

To illustrate this method, we consider the simple example%
\begin{equation*}
\int\limits_{1}^{x}\sin xdx.
\end{equation*}%
This gives%
\begin{equation}
\frac{d\xi }{dx}=\frac{-18\sin x+6\sin \left( 1\right) +6\left( x-1\right)
\cos x+3\left( x-1\right) ^{2}\sin \xi }{-\left( x-1\right) ^{3}\cos \xi }.
\label{example}
\end{equation}%
We will consider an upper limit of $x=10$ here. To find an initial value, we
choose $x_{0}\in \left( 1,10\right) ,$ such that $x_{0}$ is not too close to 
$1$ (to avoid potential problems in the denominator in (\ref{example})). We
substitute $x_{0}$ into (\ref{problem}) to get%
\begin{equation*}
\frac{\left( x_{0}-1\right) }{2}\left( \sin \left( 1\right) +\sin
x_{0}\right) +\frac{\left( x_{0}-1\right) ^{3}}{12}\sin \left( \xi \left(
x_{0}\right) \right) -\int\limits_{1}^{x_{0}}\sin xdx=0
\end{equation*}%
which we solve for $\xi \left( x_{0}\right) $ using a numerical method such
as Newton's Method or the Bisection Method. Of course, we will have to
invest some effort in obtaining an accurate value for $\int_{1}^{x_{0}}\sin
xdx,$ which can be achieved using the \textit{composite} Trapezium Rule, for
example. This effort will be well rewarded.

Once $\xi _{0}=\xi \left( x_{0}\right) $ has been determined, we solve (\ref%
{example}) subject to the initial value $\left( x_{0},\xi _{0}\right) ,$
from $x_{0}$ up to $10,$ and then from $x_{0}$ down to $1,$ using a suitable
RK method. We use a seventh-order method (RK7) \cite{butcher} in this work.
We provide an insight into the reverse RK process from $x_{0}$ down to $1$
in the Appendix.

We choose $x_{0}=5.$ We find, using composite Trapezium quadrature \cite%
{kincaid}, $\int_{1}^{5}\sin xdx=0.256640120404911,$ accurate to $\sim
10^{-15}.$ Hence, we find $\xi _{0}=3.049296665128674$ using the Bisection
Method \cite{burden}, also to an accuracy of $\sim 10^{-15}.$ Using RK7 to
solve the DE gives the results shown in the figures. In Figure (a), we show
the true value of the integral and the Trapezium approximation (dotted
line), vs $x$. Clearly, the Trapezium Rule is inaccurate. In Figure (b), we
show the error in the Trapezium Rule vs $x,$ where we see that the error is,
in places, larger than the integral itself. In Figure (c), we show $\xi
\left( x\right) .$ In Figure (d), we show the quantity%
\begin{equation*}
\left\vert \int\limits_{1}^{x}\sin xdx-\left( \underset{\text{Trapezium Rule}%
}{\underbrace{\frac{\left( x-1\right) }{2}\left( \sin \left( 1\right) +\sin
x\right) }}+\underset{\text{Error term}}{\underbrace{\frac{\left( x-1\right)
^{3}}{12}\sin \left( \xi \left( x\right) \right) }}\right) \right\vert .
\end{equation*}%
It is clear that adding the error term to the Trapezium approximation yields
values that are considerably more accurate. In a sense, we have used the
error term as a \textit{correction} term.

Considering the more exotic integrand $f\left( x\right) =x^{2}\left( \sin
x\ln \left( 2+x\right) -100x\right) ,$ which is not analytically integrable,
and using the same limits and the same $x_{0}$ as above, we find the error
in the Trapezium Rule is as high as $\sim 2\times 10^{5}.$ However, when we
add the correction term the error is reduced to $\sim 10^{-10},$ a
remarkable improvement indeed.

\section{Conclusion}

We have shown how the error term for the Trapezium Rule can be estimated, by
solving an appropriate initial value problem using a Runge-Kutta method.
Adding the error term to the Trapezium approximation yields significantly
more accurate results. Naturally, the accuracy that is attained is dependent
on the accuracy of the Runge-Kutta solution. In our example, we used a small
enough stepsize so that we achieved a high degree of accuracy. Ordinarily,
we would not know the required stepsize beforehand, and we would implement
the Runge-Kutta method with some form of error control - although not doing
that here does not detract from our central result.

\section{Appendix A}

The easiest way to implement a RK method%
\begin{align}
y_{i+1}& =y_{i}+h_{i+1}F\left( x_{i},y_{i},h_{i+1}\right)  \notag \\
h_{i+1}& =x_{i+1}-x_{i}  \label{h}
\end{align}%
in the negative $x$ direction, starting at a point $x_{0},$ is simply to
label the nodes according to%
\begin{equation*}
\left\{ \ldots ,x_{4},x_{3},x_{2},x_{1},x_{0}\right\}
\end{equation*}%
where we have $x_{0}>x_{1}>x_{2}>x_{3}>x_{4}>\ldots .$ Note that, in (\ref{h}%
), $h_{i+1}<0$ since $x_{i+1}<x_{i}.$

\section{Appendix B}

The risk of a singularity in (\ref{DE}), due to the denominator in the RHS,
can be mitigated. Since we know $f^{\prime \prime \prime }\left( x\right) ,$
it is easy to find a constant $D$ such that $f^{\prime \prime \prime }\left(
x\right) +D$ is not close to zero anywhere on $\left[ a,x\right] $ for all
values of $x$ of interest. Define 
\begin{equation*}
g\left( x\right) \equiv f\left( x\right) +\frac{Dx^{3}}{6}.
\end{equation*}%
Hence,%
\begin{equation*}
\int\limits_{a}^{x}g\left( x\right) dx=\int\limits_{a}^{x}f\left( x\right)
dx+\int\limits_{a}^{x}\frac{Dx^{3}}{6}dx
\end{equation*}%
and the Trapezium Rule applied to $g\left( x\right) $ is%
\begin{equation*}
\frac{\left( x-a\right) }{2}\left( g\left( a\right) +g\left( x\right)
\right) =\frac{\left( x-a\right) }{2}\left( f\left( a\right) +f\left(
x\right) \right) +\frac{\left( x-a\right) }{2}\left( \frac{Da^{3}}{6}+\frac{%
Dx^{3}}{6}\right) .
\end{equation*}%
This gives%
\begin{align*}
\int\limits_{a}^{x}g\left( x\right) dx-\frac{\left( x-a\right) }{2}\left(
g\left( a\right) +g\left( x\right) \right) =\text{ }&
\int\limits_{a}^{x}f\left( x\right) dx-\frac{\left( x-a\right) }{2}\left(
f\left( a\right) +f\left( x\right) \right) \\
& +\int\limits_{a}^{x}\frac{Dx^{3}}{6}dx-\frac{\left( x-a\right) }{2}\left( 
\frac{Da^{3}}{6}+\frac{Dx^{3}}{6}\right)
\end{align*}%
so that%
\begin{align*}
\int\limits_{a}^{x}f\left( x\right) dx-\frac{\left( x-a\right) }{2}\left(
f\left( a\right) +f\left( x\right) \right) =\text{ }& \left(
\int\limits_{a}^{x}g\left( x\right) dx-\frac{\left( x-a\right) }{2}\left(
g\left( a\right) +g\left( x\right) \right) \right) \\
& -\left( \int\limits_{a}^{x}\frac{Dx^{3}}{6}dx-\frac{\left( x-a\right) }{2}%
\left( \frac{Da^{3}}{6}+\frac{Dx^{3}}{6}\right) \right) .
\end{align*}%
The first term in parentheses on the RHS is the error curve that we obtain
when we apply our algorithm to $\int_{a}^{x}g\left( x\right) dx,$ and the
second term in parentheses is easy to compute exactly. The net result is the
error curve corresponding to $f\left( x\right) ,$ as desired.

For the example considered earlier, it was fortunately not necessary to use
this technique; however, had it been required, $D=2$ would have been a
suitable choice (since $f^{\prime \prime \prime }\left( x\right) =-\cos x$
in that example).

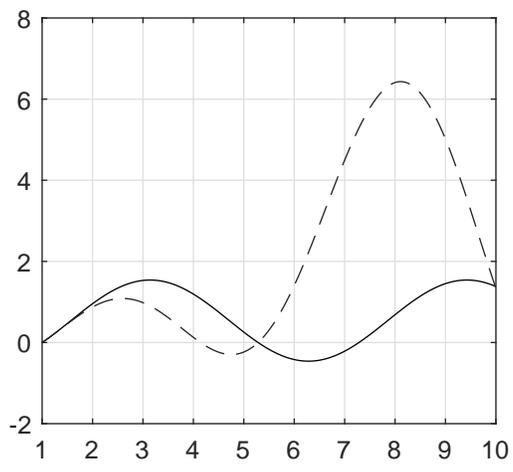
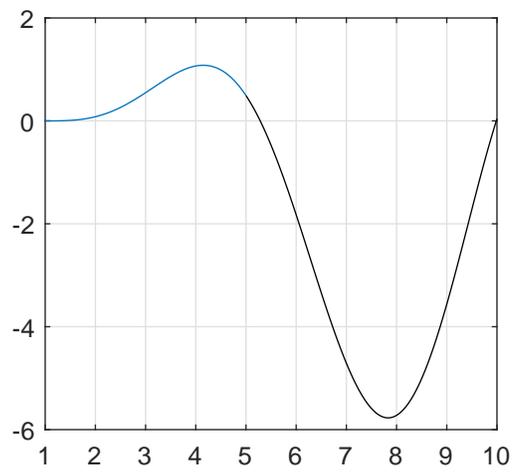
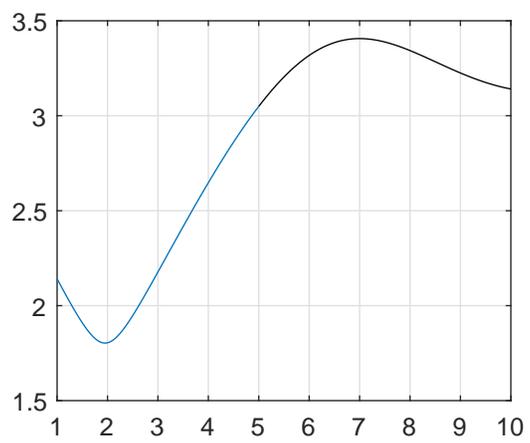
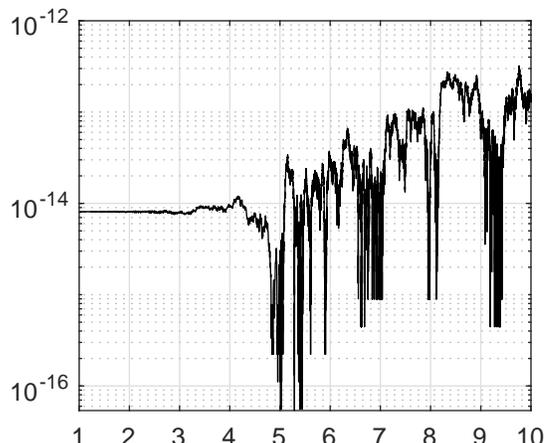
\end{document}